\documentclass{scrartcl}

\KOMAoptions{paper=a4}
\KOMAoptions{fontsize=12pt}
\KOMAoptions{draft=true}
\KOMAoptions{DIV=calc} % Berechnung des Satzspiegels auf Basis der Schriftart, Schriftgröße

\usepackage[utf8]{inputenc} % erlaubt Umlaute direkt einzugeben
\usepackage[T1]{fontenc} % erlaubt Umlaute direkt einzugeben
\usepackage[ngerman,english]{babel} % für deutsche Silbentrennung, deutsche Überschriften (z.B. Literaturverzeichnis), Datum,...
\usepackage{csquotes} % benutze \enquote{} für Anführungszeichen je nach Sprache
\usepackage{hyphenat} % advanced hyphenation: Use \hyp{} instead of -   and \fshyp{} instead of /

\usepackage{kpfonts}

\usepackage{cabin}

\usepackage{microtype} % Slightly tweak font spacing for aesthetics

\addtokomafont{title}{\rmfamily}

\usepackage{enumitem}
\setlist[enumerate]{label*=(\alph*),ref=(\alph*),itemsep=0pt,topsep=5pt}
\setlist[itemize]{itemsep=0pt,topsep=5pt}

\KOMAoptions{DIV=calc}

\usepackage{booktabs} % advanced table layout

\usepackage{graphicx} % Einbinden von Grafiken
\graphicspath{{pic/}} % für gesondertes Grafikverzeichnis

\usepackage[format=plain,labelfont={sf,footnotesize},textfont=small,margin=12pt]{caption} % Layout der Float-Unterschriften
\DeclareCaptionLabelSeparator{MySpace}{\enskip } % use half length of m as label seperator
\usepackage[style=alphabetic,
						sorting=nyt,
						maxnames=4,
						backend=biber,
						doi=true,
						isbn=false,
						clearlang=true]
						{biblatex}

\bibliography{C:/Users/jorayoga/Dropbox/Bibliography/biblio} %Neuer Laptop

\DeclareFieldFormat*{title}{\textit{#1}} % Title in italic
\renewbibmacro{in:}{} % Remove "In:" before Journal name
\DeclareFieldFormat{journaltitle}{\iffieldequalstr{journaltitle}{ArXiv e-prints}{Preprint}{#1}} % Journal name in standard font, replace ArXiv e-prints by Preprint
\DeclareFieldFormat{booktitle}{#1} % Book name in standard font
\DeclareRedundantLanguages{english,English}{english,german,ngerman,french}
\DeclareSourcemap{ 
  \maps[datatype=bibtex]{
    \map{ % removes Eprintclass, e.g. [math-AG]
		  \step[fieldset=eprintclass,null]
		}
		\map{ % only print DOI or URL if no volume information existent
      \step[fieldsource=volume,final]
      \step[fieldset=doi,null]
			\step[fieldset=url,null]
    }  
  }
}

\usepackage{hyperref} % Hyperlinks einfügen
\usepackage{xcolor} % Einfärben der Links
\hypersetup{breaklinks=true} 
\definecolor{darkblue}{RGB}{0,0,170}
\definecolor{darkred}{RGB}{200,0,0}
\hypersetup{citecolor=darkblue, urlcolor=darkblue, linkcolor=darkblue, colorlinks=true} % keine Kästen, sondern farbige Schrift
\usepackage[all]{hypcap} % Float-Links point to the top, not to the bottom, of the float

\usepackage{aliascnt} % Alias counters (in order to make autoref work)
\usepackage[amsmath,thmmarks]{ntheorem} % Theorem-Umgebungen

\newaliascnt{theorem}{theoremT}
\newtheorem {theorem}[theorem]{Theorem}
\aliascntresetthe{theorem}

\newaliascnt{proposition}{theorem}

\aliascntresetthe{proposition}

\newaliascnt{lemma}{theorem}
\newtheorem {lemma}[lemma]{Lemma}
\aliascntresetthe{lemma}

\newaliascnt{corollary}{theorem}
\newtheorem {corollary}[corollary]{Corollary}
\aliascntresetthe{corollary}

\newaliascnt{conjecture}{theorem}

\aliascntresetthe{conjecture}

\theorembodyfont{\upshape}
											 
\newaliascnt{definition}{theorem}

\aliascntresetthe{definition}

\newaliascnt{example}{theorem}
\newtheorem {example}[example]{Example}
\aliascntresetthe{example}

\newaliascnt{exercise}{theorem}

\aliascntresetthe{exercise}

\newaliascnt{goal}{theorem}

\aliascntresetthe{goal}

\newaliascnt{construction}{theorem}

\aliascntresetthe{construction}

\theoremheaderfont{\itshape}
\theorembodyfont{\upshape}

\newaliascnt{remark}{theorem}
\newtheorem {remark}[remark]{Remark}
\aliascntresetthe{remark}

\newaliascnt{convention}{theorem}

\aliascntresetthe{convention}

\newaliascnt{notation}{theorem}

\aliascntresetthe{notation}

\theoremstyle {nonumberplain}
\theoremheaderfont{\itshape}
\theorembodyfont{\upshape}
\theoremsymbol{\ensuremath{_\blacksquare}}

\newtheorem {proof}{Proof}

\usepackage {amsmath} % Mathe-Umgebungen align, gather, split, multline,...
\usepackage{tikz-cd}
\usetikzlibrary{babel}

\newcommand{\C}{{\mathbf C}}

\newcommand{\N}{{\mathbf N}}

\newcommand{\R}{{\mathbf R}}

\newcommand{\Z}{{\mathbf Z}}

\newcommand{\CP}{\mathbf{CP}}

\newcommand{\Ctor}{(\C^\times)^n}
				 
\renewcommand{\AA}{{\mathcal A}}
\newcommand{\DD}{{\mathcal D}}

\newcommand{\MM}{{\mathcal M}}

\newcommand{\Gambar}{\overline{\Gamma}}

\DeclareMathOperator{\Log}{Log}

\DeclareMathOperator{\ord}{ord}

\renewcommand{\H}{\operatorname{H}}

\begin {document}

\title {Spines for amoebas of rational curves}
\author {Grigory Mikhalkin, Johannes Rau}
\date{} 
\maketitle

\begin{abstract}
	\noindent
  To every rational complex curve $C \subset (\C^\times)^n$
	we associate a rational tropical curve $\Gamma \subset \R^n$
	so that the amoeba $\AA(C) \subset \R^n$ of $C$
	is within a bounded distance from $\Gamma$.
	In accordance with the terminology introduced in
	\cite{PR-AmoebasMongeAmpere}, we call $\Gamma$ the \emph{spine} of $\AA(C)$. 
	We use spines to describe tropical limits of sequences of rational complex curves.
\end{abstract}

{\footnote{{Research is supported in part by  the SNSF-grants 
178828, 182111  and NCCR SwissMAP}}}
{\footnote{{MSC: Primary 14H50, 14T05, 30F15}}}

\section{Introduction}

As suggested by Gelfand, Kapranov and Zelevinsky \cite{GKZ-DiscriminantsResultantsMultidimensional},
an algebraic variety $V$ in the complex torus $\Ctor=(\C\setminus\{0\})^n$
can be visualized through its {\em amoeba}.
Namely, consider the map $\Log:\Ctor\to\R^n$ defined 
by $\Log(z_1,\dots,z_n)=(\log|z_1|,\dots,\log|z_n|)$.
The image $\Log(V)$ is called the amoeba of $V$. 
It possesses many geometric properties reflecting those of $V$.
Furthermore, amoebas can be used as intermediate geometric objects
between complex and tropical varieties, cf. \cite{Mik-AmoebasAlgebraicVarietiesa}.
Passare and Rullgård \cite{PR-AmoebasMongeAmpere}
have identified a tropical variety (called {\em the Passare-Rullgård spine}) inside $\Log(V)$
in the case when $V\subset\Ctor$ is a hypersurface, i.e. $\dim V=n-1$. 

In the paper we focus on the case when $V$ is a rational curve.
In this case we associate to $V$ a tropical rational curve in $\R^n$, called {\em spine},
whose distance to $\Log(V)$ (in Hausdorff metric on sets in $\R^n$)
is universally bounded in terms of the degree of $V$. 
Our spine is not necessarily contained in $\Log(V)$.

In the case $n=2$ a rational curve $V\subset (\C^\times)^2$ is a hypersurface,
so the Passare-Rullg\aa rd spine of $V$ is also defined as a tropical curve in $\R^2$.
Nevertheless this tropical curve does not have to be a rational curve (see Remark \ref{PR-rem}).

We are freely using some basic notions from tropical geometry here. 
For details, we refer the reader to
\cite{MR-TropicalGeometry,MS-IntroductionTropicalGeometry}.

\section{The main statements}

A \emph{complex rational curve} in $\Ctor$ is a holomorphic map $f : S \to \Ctor$ from 
a Riemann sphere with $k$ labelled punctures $S = \CP^1 \setminus \{\alpha_1, \dots, \alpha_k\}$ 
to $\Ctor$.
To each puncture $\alpha_i$ we associate the integer vector $\delta(\alpha_i) \in \Z^n$ 
whose $j$-th coordinate is given by the order of vanishing 
$-\ord_{\alpha_i}(z_j \circ f)$.
The sequence of vectors $\Delta(f) = (\delta(\alpha_1), \dots, \delta(\alpha_k)) \in \Z^{k \times n}$
is called the \emph{toric degree} of $f : S \to \Ctor$.
Note that $\sum_{i=1}^k \delta(\alpha_i) = 0$.

A \emph{tropical rational curve} in $\R^n$ is a tropical morphism $h : \Gamma \to \R^n$,
where $\Gambar$ is a compact smooth rational tropical curve with $k$ labelled ends $a_1, \dots, a_k$
and $\Gamma = \Gambar \setminus \{a_1, \dots, a_k\}$.
This amounts to the following list of properties:
\begin{itemize}
	\item the graph $\Gambar$ is a tree with $k$ labelled ends $a_1, \dots, a_k$;
	\item the open subset $\Gamma$ carries a complete inner metric such that each leaf and bounded edge
	      is isometric to $[0,\infty)$ and $[0,l(e)]$, respectively. In the second case,
				$l(e) \in \R_>$ is called the \emph{length} of $e$;
	\item the map $h$ is affine on each edge;
	\item for each oriented edge $e$, the vector of derivatives $\partial h (e)$ 
	      with respect to travelling along $e$ with unit speed is integer, $\partial h (e) \in \Z^n$;
	\item at each vertex $v \in \Gamma^\circ$, if $e_1, \dots, e_k$ denote the adjacent outgoing edges,
	      the \emph{balancing condition} 
				\[
				  \sum_{i=1}^m \partial h (e_i) = 0
				\]
				is satisfied. 
\end{itemize}
Let $l_1, \dots, l_k$ denote the leaves adjacent to the ends $a_1, \dots, a_k$, oriented towards
the ends, and set $\delta(a_i) := \partial h (l_i)$, $i = 1,\dots, k$. 
The sequence of vectors $\Delta(h) = (\delta(a_1), \dots, \delta(a_k)) \in \Z^{k \times n}$
is called the \emph{toric degree} of $h : \Gamma \to \R^n$.
The balancing condition implies $\sum_{i=1}^k \delta(a_i) = 0$.

We consider the coordinate-wise logarithm map
\begin{align*} %\label{eq:} %\nonumber
  \Log \colon \Ctor &\to \R^n, \\
	       (z_1, \dots, z_n) &\mapsto (\log |z_1|, \dots, \log |z_n|). 
\end{align*}
The image of a complex curve $X$ under this map is called the \emph{amoeba} of $X$. 

\paragraph{Tropical spines} 

Our first main theorem states that the amoeba of a complex rational curve of given toric degree
can be approximated by a tropical rational curve of the same degree up to a constant which only
depends $\Delta$, but \emph{not} on the specific curve. 

Let us fix a collection of integer vectors $\Delta = (\delta_1, \dots, \delta_k)$, $\delta_i \in \Z^n$
such that $\sum_{i=1}^k \delta_i = 0$, called a \emph{toric degree} in the following.  

\begin{theorem} \label{thm:main}
  For any toric degree $\Delta$, there exists a positive constant $\epsilon = \epsilon(\Delta) \geq 0$
  having the following property. For any complex rational curve $f : S \to \Ctor$ of toric degree $\Delta$,
	there exists a tropical rational curve $h : \Gamma \to \R^n$ of toric degree $\Delta$ such that 
	\begin{align} %\label{eq:} %\nonumber
	  \Log(f(S)) \subset U_\epsilon(h(\Gamma)) & & \text{ and } & & h(\Gamma) \subset U_\epsilon(\Log(f(S))).
	\end{align}
	Here, $U_\epsilon(X)$ denotes the $\epsilon$-neighbourhood of a set $X$ in $\R^n$.
\end{theorem} 

\begin{remark} %\label{rem }
	Since all norms on $\R^n$ are equivalent, the statement of the theorem does not depend 
	on the choice of norm.
	In practice, we will work with the maximum norm $\| \,.\, \|_\infty$. 
\end{remark}

\begin{remark} %\label{rem }
\label{PR-rem}
  In \cite{PR-AmoebasMongeAmpere}, the authors associate to any complex hypersurface 
	$V_f \subset \Ctor$ a tropical hypersurface $S_f \subset \Log(V_f) \subset \R^n$, called the \emph{spine}
	of $V_f$, and show that $S_f$ is a deformation retract of $\Log(V_f)$.
	%
	%Here, $S_f$ is constructed as the hypersurface associated to the tropical polynomial
	%which is obtained as piecewise linear lower hull of the Ronkin function $N_f : \R^n \to \R$.
	The construction overlaps with ours in the case $n=2$, i.e.\ when $C = V_f$ is a planar curve.
	However, note that in general $S_f$ can be of too large genus. 
	In particular, assuming that $C$ is rational (as considered in this paper), the spine $S_f$ 
	is not necessarily rational (i.e.\ parametrised by a tropical rational curve $h : \Gamma \to \R^2$).
	A counterexample can be constructed from a counterexample to the similar statement
	that reducibility of $C$ does not imply reducibility of $S_f$. To find such an example, we may arrange
	a generic line $L_1 \subset (\C^*)^2$ and the Cremona transform of a second line $\text{Cr}(L_2)$ such that
	the union of their amoebas forms a contractible domain in $\R^2$ while the two spines $S_1$ and $S_2$ 
	intersect transversally	(in two points). In this case, among the tropical curves contained in 
	$\Log(L_1 \cup \text{Cr}(L_2))$ and of correct degree, there is a unique reducible curve (namely $S_1 \cup S_2$)
	as well as a unique curve being a deformation retract of $\Log(L_1 \cup \text{Cr}(L_2))$ (namely the spine of 
	$L_1 \cup \text{Cr}(L_2)$). Since these two curves are not equal the claim follows. 
	As mentioned above, the example can be modified to the case of an irreducible rational curve by completing
	the picture as indicated by the dashed lines. The related question to which extent $S_f$ displays 
	the singularities of $V_f$ has been studied in \cite{Lan-GeneralisationSimpleHarnack} 
	in the case of generalized simple Harnack curves.

	\begin{figure}[]
		\begin{minipage}[r]{0.5\textwidth}
			\input{pic/RPspine.TpX}
		\end{minipage}\hfill
		\begin{minipage}[l]{0.48\textwidth}
			\caption{The union of two amoebas coming from a line $L_1$ and the Cremona transform $\text{Cr}(L_2)$
																	of a second line $L_2$. In red, we depict the union of the Passare-Rullgård spines of the individual curves. 
																	The Passare-Rullgård spine of the the union $L_1 \cup \text{Cr}(L_2)$ differs from this reducible tropical curve
																	by the green edges. It gives an irreducible tropical curve which is a deformation retract of 
																	$\Log(L_1 \cup \text{Cr}(L_2))$.}
			\label{RPspine}
		\end{minipage}		
	\end{figure}%

  Despite this behaviour of Passare-Rullgård spines, one may proceed in spirit of \autoref{thm:main} 
	and try to find universal bounds $\epsilon = \epsilon(\text{NP}(f))$, only depending on the Newton polytope of
	$f$, such that 
	\begin{align} \label{eq:universalBounds} %\nonumber
	  \Log(V_f) \subset U_\epsilon(S_f) & & \text{ and } & & S_f \subset U_\epsilon(\Log(V_f)).
	\end{align}
	To our knowledge such bounds are currently not known. 
	If instead of the Passare-Rullgård spine
	the naive tropicalization of $f$ (replacing all coefficients $a_i$ by $\log|a_i|$) is used, such bounds
	have been established (at least for the first inclusion of \autoref{eq:universalBounds}) 
	in \cite{Mik-EnumerativeTropicalAlgebraic, EPR-TropicalVarietiesExponential, For-MultivariateFujiwaraBound}.
\end{remark}

\paragraph{Tropical limits}

Using \autoref{thm:main} we can describe all possible \emph{tropical limits} 
of families of rational complex curves of toric degree $\Delta$. 
Such a description is important in the context of 
correspondence theorems between complex and tropical curves.

Let $\Delta = (\delta_1, \dots, \delta_k)$ be a toric degree. 
Let $D$ be a tree with $k$ labelled leaves. We can uniquely decorate
the oriented edges of $D$ with integer vectors $\delta(e)$ 
such that
\begin{itemize}
	\item the leaf labelled by $i$ (oriented outwards) is decorated by $\delta_i$,
	\item an oppositely oriented edge $-e$ carries the vector $\delta(-e) = - \delta(e)$,
	\item around each vertex $v$, the vectors $\delta(e)$ of adjacent edges, oriented outwards,
	      sum up to zero and hence form a toric degree denoted $\Delta_v$.
\end{itemize}
A subset of vertices $S$ is called \emph{allowable} if 
there exists an assignment of non-negative non-all-zero numbers $(a(e) : e \text{ non-leaf})$
such that for any $v,w \in S$ we have 
\[
  \sum_{e \subset [v, w]} a(e) \delta(e) = 0.
\]
Here, $[v,w]$ denotes the oriented simple path from $v$ to $w$. 
A collection of toric degrees obtained as $(\Delta_v)_{v \in S}$ for an allowable vertex set $S$
is called a \emph{degeneration} of $\Delta$. An example is given in \autoref{HirzebruchCurves}.

Let $(t_m)_{k \in \N}$ be a sequence of positive real numbers converging to $+\infty$. 
Let $f_m : S_m \to \Ctor$ be a sequence of complex rational curves of 
fixed toric degree $\Delta$. We set
\[
  \AA_m := \Log_{t_m}(f_m(S_m)) = \frac{1}{\log(t_m)} \Log(f_m(S_m)).
\]
Our result describes the possible limits of such sets in the Hausdorff sense.
For precise definitions, we refer to \autoref{TropicalLimits}. 

\begin{theorem} \label{thm:TropicalLimits} \leavevmode
  \begin{enumerate}
		\item Any sequence of complex rational curves $f_m : S_m \to \Ctor$ of toric degree $\Delta$
					contains a subsequence such that the sets $\AA_m$ converge to a Hausdorff limit $A \subset \R^n$ 
					(including $A = \emptyset$).  
		\item In this case, the Hausdorff limit $A$ is of the form
					\[
						A = h_1(\Gamma_1) \cup \dots \cup h_s(\Gamma_s)
					\]
					for tropical rational curves $h_i : \Gamma_i \to \R^n$ of toric degree $\Delta_i$ 
					such that $(\Delta_1, \dots, \Delta_s)$ is a degeneration of $\Delta$. 
		\item If $s > 1$, the tropical curves $(h_1, \dots, h_s)$ can be chosen from a sublocus
					of dimension strictly less than  $n + k - 3$ in the parameter spaces of all tuples of curves 
					of degree $(\Delta_1, \dots, \Delta_s)$. 
	\end{enumerate}
\end{theorem}

Note that $n+k-3$ is the dimension of the parameter space of rational complex/tropical curves
of toric degree $\Delta$. 

The toric degree $\Delta$ also defines a homology class $\Delta_X\in H_2(X)$
in any compact toric variety $X$. This class can be obtained as the homology class
of the closure of a complex curve of toric degree $\Delta$ in $(\C^\times)^n\subset X$.
The group $H_2(X)$ also contains elements representable by curves contained in the
toric boundary $\partial X=X\setminus (\C^\times)^n$. The element $\Delta_X \in H_2(X)$ can be
represented by reducible (stable) rational curves with some components in $\partial X$.
\autoref{thm:TropicalLimits} can be used to produce examples of tropical curves that cannot
appear as tropical limits of complex curves of degree $\Delta_X$ without components in $\partial X$.

\begin{example} %\label{ex }
	Conisder the second Hirzebruch surface $\Sigma_2$ and the class $2E + 4F \in H_2(\Sigma_2)$,
	where $E$ and $F$ denote the class of the	$-2$-curve and a fibre, respectively. 
	The discriminant $\DD \subset |2E + 4F|$ consists
	of two components: The closure of the locus of irreducible rational curves,
	and the locus of reducible curves $E + |E + 4F|$. Both components have dimension $7$. 
	Tropical curves of degree $E + |E + 4F|$ 
	(actually, since we restrict to $\R^2$, in $|E + 4F|$)
	also form a $7$-dimensional family. By \autoref{thm:TropicalLimits},
	the ones that appear as limits of irreducible complex curves form 
	a subfamily of dimension at most $6$. 
	\autoref{HirzebruchCurves} shows examples of curves 
	which appear or do not appear as such limits. 
\end{example}

% TpX Bild einfügen
	\begin{figure}[ht]
		\centering
		\input{pic/HirzebruchCurves.TpX}
		\caption{Three tropical curves of type $|E + 4F|$ in $\Sigma_2$. 
	Since $((0,-1)^4, (1,2), (0,1)^2,(-1,0))$ is not a degeneration of 
	$\Delta = ((0,-1)^4, (1,2)^2, (-1,0)^2)$, $C_1$ cannot occur as limit of 
	(irreducible) complex rational curves in $|2E + 4F|$. 
	Conversely, $C_2$ and $C_3$ occur as limits of tropical curves 
	of toric degree $\Delta$ as shown on their right hand side. 
	The degree of $C_2$ is $((0,-1)^4, (1,2), (0,2),(-1,0))$ and may be 
	written as a degeneration of $\Delta$ using a single allowable vertex $S$. 
	The pair of degrees of the reducible curve $C_3$ is 
	$((0,1), (0,-1))$ and $((0,-1)^3, (1,2), (0,1),(-1,0))$.
	It is a degeneration of $\Delta$ induced by the tree $D$ and
	allowable vertices $S$ (in black) displayed on the bottom left hand side.
	The path connecting the two vertices in $S$ consists of two edges of
	direction $(0,\pm 1)$, and we can choose $a(e) = 1$ for these two edges and
	$a(e) = 0$ for all others. 
	}
		\label{HirzebruchCurves}
	\end{figure}%

\section{Spines of lines}

For every $n \in \N$, we set
\[
  \Delta_n := (-e_0, -e_1, \dots, -e_n),
\]
where $e_1, \dots, e_n$ denotes the standard basis of $\R^n$ and $-e_0 = e_1 + \dots e_n = (1, \dots, 1)$. 
Complex and tropical curves of degree $\Delta_n$ are called \emph{(non-degenerate) lines}. 
For lines, we number the punctures, ends, and leaves, respectively, from $0$ to $n$.

\paragraph{Complex lines} By choosing a coordinate $z$ for $\CP^1$ such that $\alpha_0= \infty$, 
we can parametrize any complex line 
$L \subset \Ctor$ by a map
\begin{align} %\label{eq:} %\nonumber
  f \colon \CP^1 \setminus \{\infty, \alpha_1, \dots, \alpha_n\} &\to \Ctor, &
	z &\mapsto (\kappa_1 (z - \alpha_1), \dots, \kappa_n (z - \alpha_n)).
\end{align}
We call $L$ \emph{calibrated} if $\kappa_1 = \dots = \kappa_n$. 

\paragraph{Tropical lines} \label{ItemsTropicalLines}
Let us recall the basic properties of tropical lines:
\begin{itemize}
	\item If $h : \Gamma \to \R^n$ is a tropical line, then $h$ is injective. 
	      Indeed, the balancing condition implies that, as we follow the 
				path from $a_i$ to $a_0$ with unit speed, the function $x_i \circ h$
				has constant derivative $1$. Since any point $p \in \Gamma$ lies 
				on at least one such path, the injectivity follows.
	\item Throughout the following, we will identify $\Gamma$ with its image 
				and use the notation $\Gamma \subset \R^n$ (suppressing $h$).
	\item Given $p = (p_1, \dots, p_n) \in \Gamma$, let $e$ denote the oriented edge pointing from $p$ towards
	      $a_0$. Then the direction vector 
				$\partial h (e) \in \Z^n$ has only $0$ and $1$ as entries. 
				A coordinate $x_i$ corresponding to an entry $1$ is called 
				a \emph{local coordinate} for $\Gamma$ at $p$.
				Given $k \in \R$, the linear tropical polynomials 
				$\mu(x) = k + \max\{x_i-p_i, 0\}$,
				for any local coordinate $x_i$, 
				restrict to the same function on $\Gamma$. 
				The \emph{linear modification} of $\Gamma$ at $p$ (of height $k$) is
				the unique line $\widetilde{\Gamma} \subset \R^{n+1}$ which contains
				the graph of $\mu$. 
				More concretely, $\widetilde{\Gamma}$ is the union 
				of the graph of $\mu$ with the
				ray in direction $-e_{n+1}$ 
				emanating from $(p,k)$. 
	\item The inverse operation to modification is called contraction. 
				Let $\pi : \R^n \to \R^{n-1}$ be the projection forgetting $x_n$.
	      Then the image $\Gamma' = \pi(\Gamma)$ of any line $\Gamma \subset \R^n$ 
				is a line in $\R^{n-1}$, called the \emph{contraction} of $\Gamma$ (along $x_n$).
				Let $p \in \Gamma'$ be the image of the contracted leaf $l_n$. 
				Then $\Gamma$ is the linear modification of $\Gamma'$ at $p$ (for a suitable height $k$). 
				In particular, the contraction map $\pi : \Gamma \to \Gamma'$ is a bijection when 
				restricted to $\Gamma \setminus l_n^\circ$.
	\item A tropical line $\Gamma$ is \emph{calibrated} if the leaf $l_0$
	      is contained in the (usual) line $\R e_0$ (emanating from the origin). 
				Given $p \in \Gamma$, note that $\Gamma$ is calibrated if and only if 
				$p_i = p_j$ for any two local coordinates	$x_i$ and $x_j$ at $p$. 
				Moreover, the modification of a calibrated
				line is calibrated if and only, in the notation from above, $k=p_i$
				and hence $\mu(x) = \max\{x_i, p_i\}$.
\end{itemize}

\paragraph{Spines for amoebas of lines} We consider the (shifted) geometric series 
\[
  \epsilon_n = 2 \log(2) \sum_{i=0}^{n-2} 3^i = \log(2) (3^{n-1}-1)
\]
with initial value $\epsilon_1 = 0$.
Note that $\epsilon_n = 3 \epsilon_{n-1} + 2\log(2)$ for all $n \in \N$.

\begin{theorem} \label{thm:mainlines}
  Let $L \subset \Ctor$ be a complex line. Then there exists a tropical line $\Gamma \subset \R^n$ and 
	a map $\phi : L \to \Gamma$ such that 
	\[
	  \| \Log(q) - \phi(q) \|_\infty \leq \epsilon_n
	\]
	for all $q \in L$. Moreover, if $L$ is calibrated, there exists a calibrated $\Gamma$
	such that the statement holds.
\end{theorem}

\begin{figure}
	\begin{minipage}[c]{0.52\textwidth}
		\includegraphics[width=\textwidth]{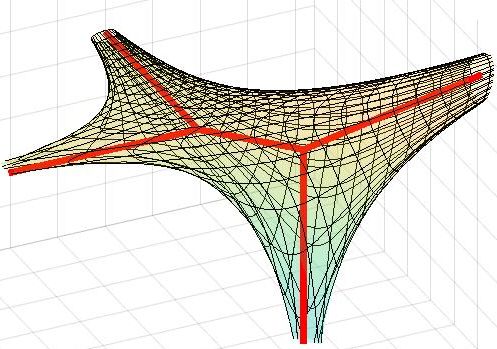}
	\end{minipage}\hfill
	\begin{minipage}[r]{0.45\textwidth}
		\caption{The amoeba of a complex line $L \subset (\C^\times)^3$ together with 
						 an approximating tropical line $\Gamma \in \R^3$. The line $L$ is 
						parametrised by	$z \mapsto (z, z+1, z - 2i)$. The vertices of $\Gamma$
						are $(0,0,0)$ and $(\log(2), \log(2),0)$.}
		\label{amoebaline}
	\end{minipage}		
\end{figure}%

\begin{proof}
  We prove the statement for calibrated lines by induction on $n$. 
	The general statement obviously follows from the calibrated case 
	after applying translations in $\Ctor$ and $\R^n$. 
	
	For $n=1$, we have $L = \C^\times$ and hence
	$\Gamma = \R$ and $\phi = \Log : \C^\times \to \R$ satisfy the requirements. 
	
	For the induction step $n-1 \to n$, let us start with a given calibrated complex line $L \subset (\C^\times)^{n}$. 
	We denote by $L' \subset (\C^\times)^{n-1}$ the calibrated complex line obtained as the closure 
	of the image of $L$ under the projection forgetting the last coordinate $z_{n}$. 
	The closure contains the point $w = (w_1, \dots, w_{n-1})$ corresponding to
	the puncture $\alpha_n$.
	Since $L$ is calibrated, the coordinates on $L$ are related by $z_n = z_i - w_i$ for $i = 1, \dots, n-1$.
	
	By the induction assumption, there
	exists a calibrated tropical line $\Gamma' \subset \R^{n-1}$ and a map $\phi' : L' \to \Gamma'$ 
	such that $\| \Log'(q) - \phi'(q) \|_\infty \leq \epsilon_{n-1}$ for all $q \in L'$.
	Here, we use $\Log$ and $\Log'$ to denote the log map on $n$ and $n-1$ variables, respectively. 
	We set $p = (p_1, \dots, p_{n-1}) := \phi'(w)$.
	We define the tropical line $\Gamma \subset \R^{n}$ as the modification 
	of $\Gamma'$ at $p$ corresponding to the function $\mu(x) = \max\{x_i, p_i\}$
	for a local coordinate $x_i$ at $p$. 
	By the remarks on page \pageref{ItemsTropicalLines}, $\Gamma$ is calibrated and does not depend
	on the choice of local coordinate $x_i$. 
	For now, let us fix such $x_i$. 
	
	In the next step, we define the map $\phi : L \to \Gamma$. We distinguish two cases depending on whether 
	$\phi'(q)$ is close to $p$ or not. For $q = (q_1, \dots, q_n) \in L$ we set
	\[
	  \phi(q) := 
		\begin{cases}
		  \big(p,\min\{\log|q_n|, p_i\}\big) & \text{if } \|\phi'(q) - p \|_\infty \leq 2 \epsilon_{n-1} + \log(2), \\
			\big(\phi'(q), \mu(\phi'(q))\big)  & \text{otherwise}.
		\end{cases}
	\]
	Note that $\phi(q) \in \Gamma$ by construction. It remains to prove that 
	$\| \Log(q) - \phi(q) \|_\infty \leq \epsilon_n$ for all $q \in L$. 
	Before continuing, let us collect two consequences of the induction assumption for reference:
	\begin{align}  %\nonumber
		\big|\log|q_i| - \phi'(q)_i\big| &\leq \epsilon_{n-1}, \label{eq:indass1} \\
		\big|\log|w_i| - p_i\big| &\leq \epsilon_{n-1}. \label{eq:indass2}
	\end{align}
	We proceed in several cases.
	
	\textit{Case 1}
	Assume that $\|\phi'(q) - p \|_\infty \leq 2 \epsilon_{n-1} + \log(2)$.
	Since this implies
	\[
	  \|\Log'(q) - p\|_\infty \leq \|\Log'(q) - \phi'(q)\|_\infty + \|\phi'(q) - p\|_\infty \leq 3 \epsilon_{n-1} + \log(2) < \epsilon_n,
	\]
	by the definition of $\phi(q)$ it suffices to show $\log|q_n| \leq p_i + \epsilon_n$.
	To do so, we apply the case assumption again to the $i$-th coordinates, 
	providing
	\begin{equation} \label{eq:ass3} %\nonumber
		|\phi'(q)_i - p_i| \leq 2 \epsilon_{n-1} + \log(2).
	\end{equation}
	Combining \autoref{eq:indass1} and \autoref{eq:ass3}, we get 
	\begin{equation} \label{eq:intermediate} %\nonumber
		\big|\log|q_i| - p_i\big| \leq 3 \epsilon_{n-1} + \log(2) 
	\end{equation}
	and hence 
	\begin{equation} %\label{eq:} %\nonumber
		|q_n| = |q_i - w_i| \leq |q_i| + |w_i| \leq (2 e^{3 \epsilon_{n-1}} + e^{\epsilon_{n-1}}) e^{p_i} 
		  \leq 3 e^{3 \epsilon_{n-1}} e^{p_i} < e^{\epsilon_n} e^{p_i}.
	\end{equation}
	Here, the second inequality uses \autoref{eq:indass2} and \autoref{eq:intermediate} and the third inequality 
	follows from 
	$e^{3 \epsilon_{n-1}} \geq e^{\epsilon_{n-1}}$ since $\epsilon_{n-1} \geq 0$.
	Finally, this implies 
	\begin{align} %\label{eq:} %\nonumber
		\log|q_n| = \log|q_i - w_i| \leq p_i +\epsilon_n, 
	\end{align}
	as required.

	\textit{Case 2} 
	Let us now assume $\|\phi'(q) - p \|_\infty > 2 \epsilon_{n-1} + \log(2)$.
	By definition of $\phi(q)$, we need to show that $|\log|q_n| - \mu(\phi'(q))| \leq \epsilon_n$. 
	We subdivide this case further as follows (see \autoref{ProofCases}):
	
	% TpX Bild einfügen
	\begin{figure}[b]
		\centering
		\input{pic/ProofCases.TpX}
		\caption{The three subcases 2.1--3 in the proof of \autoref{thm:mainlines} with $\delta = 2 \epsilon_{n-1} + \log(2)$.}
		\label{ProofCases}
	\end{figure}%

	\textit{Subcase 2.1} 
	There exists $i \in \{1,\dots, n-1\}$ such that $\phi'(q)_i - p_i > 2 \epsilon_{n-1} + \log(2)$.
	Note that when following the path in $L'$ from $p$ to $\phi'(q)$, any coordinate increases at most 
	as much as the local coordinates at $p$. 
	Thus we may assume without loss of generality that $x_i$
	is a local coordinate at $p$ and hence $\mu(\phi'(q)) = \phi'(q)_i$. 
	Using Equations \ref{eq:indass1} and \ref{eq:indass2}, we obtain
	\[
	  \log|q_i| - \log|w_i| > \log(2)
	\]
	or, equivalently, $|q_i| > 2 |w_i|$. 
	The triangle inequalities for $q_n = q_i - w_i$ give
	\begin{align} %\label{eq:} %\nonumber
		|q_n| &\leq |q_i| + |w_i| < |q_i| + \frac{1}{2} |q_i| < 2 |q_i|, \nonumber \\
		|q_n| &\geq |q_i| - |w_i| > |q_i| - \frac{1}{2} |q_i| > \frac{1}{2} |q_i|, \nonumber
	\end{align}
	and hence
	\[
	  \big|\log|q_n| - \log|q_i|\big| < \log(2).
	\]
	Together with \autoref{eq:indass1}, we get $|\log|q_n| - \phi'(q)_i| < \epsilon_{n-1} + \log(2) < \epsilon_n$.
	
	\textit{Subcase 2.2} 
  There exists a local coordinate $x_i$ at $p$ such that $p_i - \phi'(q)_i  > 2 \epsilon_{n-1} + \log(2)$. 
	The reciprocal previous argument implies
	\[
	  \big|\log|q_n| - \log|w_i|\big| < \log(2)
	\]
	and hence $|\log|q_n| - p_i| < \epsilon_{n-1} + \log(2) < \epsilon_n$,
	and we are done.
	
	\textit{Subcase 2.3} 
	We have $p_i - \phi'(q)_i  > 2 \epsilon_{n-1} + \log(2)$ for 
	some $i$ and none of the previous subcases occurs.
	In this case, the subtree of $\Gamma'$ spanned by 
	$p$, $\phi'(q)$ and $l_0$ (the leaf corresponding to $e_0$) 
	contains a unique three-valent vertex 
	$r = (r_1, \dots, r_n) \in \Gamma'$. Alternatively, $r$ can be described
	as the point on the path from $p$ to $\phi'(q)$ at which the coordinate $x_i$ starts to decrease. 
	In particular, $r_i = p_i$. 
	Note that $\| r - p \|_\infty \leq 2 \epsilon_{n-1} + \log(2)$, since otherwise
	this would imply the existence of a coordinate satisfying the conditions of the first subcase.
	Let $x_j$ be a local coordinate at $p$. 
	Then $\mu(\phi'(q)) = \phi'(q)_j = r_j$ by construction of $r$. 
	Moreover, both $x_i$ and $x_j$ are local coordinates 
	at $r$. Since $\Gamma'$ is calibrated, this implies $r_i = r_j$.
	The estimate $|\log|q_n| - p_i| < \epsilon_{n-1} + \log(2)$
	from the second subcase is still valid, so we can combine these equations
	to 
	\begin{multline}
	  \big|\log|q_n| - \mu(\phi'(q))\big| = \big|\log|q_n| - r_i\big| \\
		  \leq \big|\log|q_n| - p_i\big| + \big|p_i - r_i\big| < 3 \epsilon_{n-1} + 2 \log(2) = \epsilon_n.
	\end{multline}
	This finishes the third subcase and hence completes the proof.
\end{proof}

\begin{remark} %\label{rem }
  We made no serious attempt to reach optimality of $\epsilon_n$ in any sense. For example,
	$\epsilon_2 = 2\log(2)$ can obviously be improved to $\log(2)$ (even with respect to the Euclidean metric). 
	Note also that except for the trivial case $n = 1$	the proof in fact yields the strict inequality 
	$\| \Log(q) - \phi(q) \|_\infty < \epsilon_n$.
\end{remark}

\autoref{thm:mainlines} clearly implies $\Log(L) \subset U_{\epsilon_n}(\Gamma)$. 
To prove $\Gamma \subset U_{\epsilon_n}(\Log(L))$, we upgrade the statement to 
show surjectivity of $\phi$ up to small neighbourhoods around the vertices of $\Gamma$.

\begin{theorem} \label{thm:SurjectivePhi}
  The map $\phi : L \to \Gamma$ in \autoref{thm:mainlines} can be chosen such that 
	\[
	  \Gamma \setminus \bigcup_{v \text{ vertex}} U_{\epsilon_n}(v) \subset \phi(L).
	\]
\end{theorem}

\begin{proof}
  As before, we may restrict to the calibrated case. We use the same induction as 
	in \autoref{thm:mainlines}. 
	For $n=1$, the $\phi = \Log : \C^\times \to \R$ is obviously surjective. 
	For the induction step $n-1 \to n$, we use the same notation as before
	and set $R = \Gamma \setminus \bigcup_{v} U_{\epsilon_n}(v)$
	and $R' = \Gamma' \setminus \bigcup_{v'} U_{\epsilon_{n-1}}(v')$.
	The additional induction assumption is $R' \subset (\phi'(L'))$.	
	
	Clearly $\pi(R\setminus l_n) \subset R'$. Moreover, for any $q$ 
	with $\phi'(q) \in \pi(S\setminus l_n)$, the \enquote{otherwise}-case in the 
	definition of $\phi$ is used. By the induction assumption, we conclude
	$R\setminus l_n \subset \phi(L)$. 
	
	It remains to show that a point in $l_n$ with last coordinate lower or equal
	than $p_i - \epsilon_n$ lies in $\phi(L)$. Here $x_i$ is a local coordinate for $p$. 
	In fact, we will prove the stronger statement that for any $q \in L$ with 
	$\log|q_n| \leq p_i - \epsilon_n$, the \enquote{if}-case in the definition of
	$\phi$ takes effect. First, note that $p_i \leq p_j$ for all $j$ since $\Gamma'$ is calibrated. 
	It follows that $\log|q_n| \leq \log|w_j| + \epsilon_{n-1} - \epsilon_n$ for all $j$.
	Since $\epsilon_n - \epsilon_{n-1} > \log(2)$, we get $|q_n| < |w_j| / 2$. 
	As in previous arguments, this implies $|\log|q_j| - \log|w_j|| < \log(2)$
	and hence $|\phi(q)_j - p_j| < 2 \epsilon_{n-1} + \log(2)$. This shows
	$\|\phi'(q) - p \|_\infty \leq 2 \epsilon_{n-1} + \log(2)$ and finishes the proof. 
\end{proof}

\begin{remark} \label{rem:BetterPhi}
  With little extra effort, the induction argument can be modified to 
	construct a map $\phi : L \to \Gamma$ with the following properties.
	\begin{enumerate}
		\item The map $\phi$ is continuous, proper and surjective. 
		\item For all $q \in L$ we have $\| \Log(q) - \phi(q) \|_\infty \leq \epsilon'_n$.
		\item For any $p \in \Gamma$ in the interior of an edge $e$, the preimage 
		      $\phi^{-1}(p)$ is a smoothly embedded circle in $L$ and the homology class
					$[\phi^{-1}(p)] \in \H_1(\Ctor, \Z) = \Z^n$ is equal to the direction 
					vector of $e$ (for compatible orientations of $\phi^{-1}(p)$ and $e$). 
		\item For any vertex $v \in \Gamma$, the preimage $\phi^{-1}(v) \subset L$ is a compact 
		      surface with boundary. The boundary components are in bijection (given by homology classes) 
		      with the edges adjacent	to $v$. 
	\end{enumerate}
	Here, the value $\epsilon'_n$ can be defined by the recursion 
	$\epsilon'_n = 5 \epsilon'_{n-1} + \log(5)$. 
	The induction step can then be modified as follows: 
	Choose $2 \epsilon'_{n-1} + \log(2) \leq \delta < 2 \epsilon'_{n-1} + \log(5)/2$ 
	such that $\partial U_\delta(p)$ does not contain vertices of $\Gamma'$ 
	and set 
	\[
	  \phi(q) := 
		\begin{cases}
		  \big(p,\log|q_n|\big) & \text{if } \log|q_n| \leq p_i - \delta, \\
			\big(\phi'(q), \mu(\phi'(q))\big)  & \text{if } \phi'(q) \notin U_\delta(p).
		\end{cases}
	\]
	It remains to extend $\phi$ to 
	\[
	  B = \{q \in L : \phi'(q) \in \overline{U_\delta(p)} \text{ and } \log|q_n| \geq p_i - \delta\},
	\]
	which is a connected surface with boundary in $L$ whose boundary components are in bijection with 
	$\Gamma \cap \partial U_\delta((p,p_i))$. 
	It is clear that a map $\phi_B : B \to \Gamma \cap \overline{U_\delta((p,p_i))}$ 
	satisfying properties (a), (c) and (d) exists. 
	Property (b) then follows from previous arguments and  
	\begin{align*} %\label{eq:} %\nonumber
	  \|\pi(\phi(q) - \Log(q))\|_\infty 
		  & \leq \|\pi(\phi(q)) - p\|_\infty + \|p - \phi'(q)\|_\infty + \|\phi'(q) - \Log'(q)\|_\infty \\
			& \leq \delta + \delta + \epsilon'_{n-1} < 5 \epsilon'_{n-1} + \log(5) = \epsilon'_n.
	\end{align*}
  Using $\phi_B$ to extend $\phi$ to $L$, we obtain a function which satisfies (a) -- (d).
\end{remark}

\section{Spines of rational curves}

Let $\Delta = (\delta_0, \dots, \delta_k)$ be a toric degree in dimension $n$. 
We denote by $\psi_\Delta : \R^k \to \R^n$ the linear map which sends the standard basis
vector $e_i$ to $-\delta_i$ for all $i=1,\dots,k$ (this implies $e_0 \to -\delta_0$). 
In this section, we assume that $\Delta$ is \emph{non-degenerate}, that is to say, 
the map $\psi_\Delta$ is surjective. 

Let $\Psi_\Delta : (\C^\times)^k \to \Ctor$ denote the torus homomorphism which is 
the exponential of $\phi_\Delta$ (hence also surjective). 
In other words, the diagram
\begin{equation} %\label{eq:} %\nonumber
\begin{tikzcd}
	(\C^\times)^k \arrow[r, "\Psi_\Delta"] \arrow[d, "\Log"']
	  & \Ctor \arrow[d, "\Log"] \\
	\R^k \arrow[r, "\psi_\Delta"]
	  & \R^n
\end{tikzcd}
\end{equation}
commutes. 

The following lemmas state that complex and tropical rational curves of toric degree $\Delta$ can 
be represented as images of lines under $\Psi_\Delta$ and $\psi_\Delta$, respectively. 

\begin{lemma} \label{lem:linesratlcmpl}
  Given a complex line $L \subset (\C^\times)^k$, the map 
	\[
	  f = \Psi_\Delta|_L : L \to \Ctor
	\]
	is a complex rational curve
	of toric degree $\Delta$. Any complex rational curve	of toric degree 
	$\Delta$ can be represented in such a way.  
	Two lines $L,L'$ provide the same rational curve if and only if $L = wL'$ for some $w \in \ker\Psi_\Delta$.
\end{lemma}

\begin{proof}
  The uniqueness up to $\ker\Psi_\Delta$ is obvious.
	Using coordinates $\delta_i = (\delta_i^1, \dots, \delta_i^n)$, the map $\Psi_\Delta$ is given by
	\[
	  z'_j = z_1^{-\delta_1^j} \cdots z_k^{-\delta_k^j}.
	\]
	This implies $-\ord_{\alpha_i}(z_j \circ f) = \delta_i^j$, as required. 
	
	Let $f : S \to \Ctor$ be a complex rational curve of toric degree $\Delta$. Up to isomorphism,
	we may assume $S = \CP^1 \setminus \{\infty, \alpha_1, \dots, \alpha_n\}$, with affine coordinate $z$.
	By definition of toric degree, we have
	\[
	  z_j \circ f = \kappa_j (z - \alpha_1)^{-\delta_1^j} \cdots (z - \alpha_k)^{-\delta_k^j}
	\]
	for some constant $\kappa_j \in \C^\times$.
	Pick a preimage $(\lambda_1, \dots, \lambda_k)$ of $(\kappa_1, \dots, \kappa_n)$ under
  $\Psi_\Delta$.
	Then $f$ factors through $\Psi_\Delta$ by the line
	\begin{align} %\label{eq:} %\nonumber
		S &\to (\C^\times)^k, \\
		z &\mapsto (\lambda_1 (z-\alpha_1), \dots, \lambda_k (z-\alpha_k)).
	\end{align}
\end{proof}

\begin{lemma} \label{lem:linesratltrop}
  Given a tropical line $\Gamma \subset \R^k$, the map 
	\[
	  h = \psi_\Delta|_\Gamma : \Gamma \to \R^n
	\]
	is a tropical rational curve	of toric degree $\Delta$. 
	Up to isomorphism, any tropical rational curve of toric degree $\Delta$ can be represented in such a way.  
	Two lines $\Gamma,\Gamma'$ provide the same rational curve if and only if $\Gamma = x\Gamma'$ for some $x \in \ker\psi_\Delta$.
\end{lemma}

\begin{proof}
  Let $\Gamma \subset \R^k$ be a tropical line. 
  Since $\psi_\Delta$ is linear, $\psi_\Delta|_\Gamma : \Gamma \to \R^n$ is clearly a tropical morphism.
	Moreover, the degree requirements are satisfied since $\psi_\Delta$ maps $-e_i \to \delta_i$ for
	$i = 0, \dots, k$. 
	
	Given a tree $\Gamma$ with complete inner metric and $k+1$ leaves $l_0, \dots, l_k$, 
	an arbitrary base point $p_0 \in \Gamma$ and
	toric degree $\Delta$, the set of tropical rational curves $h : \Gamma \to \R^n$ of toric degree $\Delta$ 
	is in bijection	to $\R^n$ via $h \mapsto h(p_0)$. Indeed, since $\Gamma$ is a tree and 
	since the direction vectors $\partial h (l_i)$ are
	fixed by $\Delta$, the balancing condition recursively prescribes all direction vectors $\partial h(e)$. 
	To fix $h : \Gamma \to \R^n$, it hence suffices to fix the image of a single point. 
	
	Let $h : \Gamma \to \R^n$ be a tropical rational curve of toric degree $\Delta$ with base point $p_0$. 
	Choose a point $x \in \R^k$ such that $\psi_\Delta(x) = h(p_0)$. Applying the previous discussion
	to $\Delta_k$, there exists a unique tropical line $g : \Gamma \to \R^k$ such that $g(p_0) = x$.
	Moreover, by construction we have $\psi_\Delta(\partial g(e)) = \partial h(e)$ for any edge $e$ of $\Gamma$. 
	Hence, $f = \psi_\Delta \circ g$, as required. 
	The uniqueness property also follows easily from the previous discussion.
\end{proof}

We are now ready to prove the main theorem.

\begin{proof}[\autoref{thm:main}]
  Given a toric degree $\Delta$ consisting of $k+1$ vectors, we set 
	$\epsilon' = \epsilon_{k} \cdot N(\Delta)$, where
	\[
	  N(\Delta) = \| \psi_\Delta \|_\infty = \max\left\{ \frac{\|\psi_\Delta(x)\|_\infty}{\|x\|_\infty} : 0 \neq x \in \R^k \right\}.
	\]
	Let $f : S \to \Ctor$ be a complex rational curve  of toric degree $\Delta$.
	By \autoref{lem:linesratlcmpl}, we may assume that $S \cong L \subset (\C^\times)^k$ is a complex line 
	and $f = \psi_\Delta(x)|_S$. By \autoref{thm:mainlines}, there exists a tropical line $\Gamma \subset \R^k$
	and a map $\phi : L \to \Gamma$ such that $\| \Log(q) - \phi(q) \|_\infty \leq \epsilon_k$ for all $q \in L$.
	By \autoref{lem:linesratltrop}, $h = \psi_\Delta|_\Gamma : \Gamma \to \R^n$ is a tropical rational curve
	of toric degree $\Delta$. 
	The situation can be summarized in the following diagram (whose left hand side is only commutative up to $\epsilon_k$):
	\begin{equation} %\label{eq:} %\nonumber
	\begin{tikzcd}
		L  \arrow[r, phantom, "\subset" description] \arrow[d, "\phi"'] \arrow[rr, "f" description, bend left=30, start anchor=north east, end anchor=north west]
		            \arrow[dr, phantom, "{\scriptstyle \leq\epsilon_k}" description]
		      &[-2ex] (\C^\times)^k \arrow[r, "\Psi_\Delta"'] \arrow[d, "\Log"]    
					& \Ctor \arrow[d, "\Log"] \\
		\Gamma \arrow[r, phantom, "\subset" description] \arrow[rr, "h" description, bend right=30, start anchor=south east, end anchor=south west]
		      & \R^k          \arrow[r, "\psi_\Delta"]                       
					& \R^n
	\end{tikzcd}
	\end{equation}
	Hence, for all $q \in L$,
	\begin{equation} %\label{eq:} %\nonumber
		\|\Log(f(q)) - h(\phi(q))\|_\infty = \|\psi_\Delta(\Log(q) - \phi(q))\|_\infty 
		  \leq N(\Delta) \epsilon_k = \epsilon',
	\end{equation}
	which implies $\Log(f(L)) \subset U_{\epsilon'}(h(\Gamma))$.
	
	Set $R = \Gamma \setminus \bigcup_{v} U_{\epsilon_n}(v)$.
	By \autoref{thm:SurjectivePhi}, we have 
	$h(R) \subset h(\phi(L)) \subset U_{\epsilon'}(\Log(f(L)))$.
	Finally, for $p \in \Gamma \setminus R$, there
	exists $p' \in R$ such that $\|p - p'\|_\infty < (k-1)\epsilon_k$,	
	since $\Gamma$ has $k-1$ vertices. 
	Choose $q' \in L$ with $\phi(q') = p'$.
	Then
	\[
	  \|h(p) - \Log(f(q'))\|_\infty \leq \|h(p) - h(p')\|_\infty + \|h(p') - \Log(f(q'))\|_\infty
		                              < k \epsilon'.
	\]
	Hence, for $\epsilon = \epsilon(\Delta) = k \cdot \epsilon'$ we proved
	$h(\Gamma) \subset U_{\epsilon}(\Log(f(L)))$, which finishes the proof. 
\end{proof}

\begin{remark} %\label{rem }
  Clearly, \autoref{rem:BetterPhi} can be extended to the general case in the sense that for any 
	complex rational curve $f : S \to \Ctor$ of toric degree $\Delta$,
	there exists a tropical rational curve $h : \Gamma \to \R^n$ of toric degree $\Delta$ 
	and a map $\phi : S \to \Gamma$ which satisfies properties (a) -- (d)
	(after substituting $\|\Log(f(q)) - h(\phi(q))\|_\infty \leq \epsilon'(\Delta)$ 
	and $f_*[\phi^{-1}(p)]$ at the obvious places). Here,  $\epsilon'(\Delta) = 
	N(\Delta) \epsilon'_k$.
\end{remark}

\section{Tropical limits of amoebas} \label{TropicalLimits}

Given two subsets $A,B \subset \R^n$, we set the \emph{Hausdorff distance}
of $A$ and $B$ to 
\[
  d(A,B) = \inf\{\delta : A \subset U_\delta(B), B \subset U_\delta(A)\}.
\]
Note that $d(A,B)$ can be infinite in general. 
If we restrict to non-empty closed 
subsets of a compact set $K \subset \R^n$, 
then $d(A,B) \in \R_\geq$ and the Hausdorff distance defines
a metric. 
A sequence of subsets $A_m \subset \R^n$ \emph{converges to the Hausdorff limit} $A \subset \R^n$
if $A$ is closed and for any compact set $K \subset \R^n$ the sequence $d(A_m \cap K, A \cap K)$ converges to $0$.
In this case $A = \lim A_m$ is unique, since it is unique on each compact $K$.  
Note that we include the case $A = \emptyset$, which is to say, 
for any compact $K \subset \R^n$, 
there exists $k_0 \in \N$ such that $A_m \cap K = \emptyset$ for all $k \geq k_0$. 

Let $f_m : S_m \to \Ctor$ be a sequence of rational complex curves as in the assumptions
of \autoref{thm:TropicalLimits}.  
By \autoref{thm:main}, there exists a sequence of tropical rational curves
$h_m : \Gamma_m \to \R^n$ of toric degree $\Delta$ such that 
\begin{align} %\label{eq:} %\nonumber
	\AA_m \subset U_{\epsilon/\log(t_m)}(h_m(\Gamma_m)) & & \text{ and } & & h_m(\Gamma_m) \subset U_{\epsilon/\log(t_m)}(\AA_m).
\end{align}
This implies that the sequence of Hausdorff distances $d(\AA_m, h_m(\Gamma_m))$ converges to zero. 
Obviously, this is still true after restricting to compact subsets $K$. We get the following corollary.

\begin{corollary} \label{cor:TropicalComplexLimits}
  The sequence $\AA_m$ converges to the Hausdorff limit $A$ if and only if $h_m(\Gamma_m)$
	converges to the Hausdorff limit $A$. 
\end{corollary}

In other words, \autoref{thm:main} reduces the proof of \autoref{thm:TropicalLimits} to the study of Hausdorff limits
of tropical curves. 

Fix a toric degree $\Delta = (\delta_1, \dots, \delta_k)$ in $\R^n$, 
an (abstract) tree $G$ with $m$ leaves labelled 
by $\{1, \dots, m\}$ and a marked vertex $v_0 \in G$.

In analogy to our conventions for $\Gamma$, the leaves are considered
to be half-edges without one-valent end vertices.
Let us furthermore assume that $G$ does not contain two-valent vertices except for the case
$- \bullet - $. Then the space $\MM(\Delta, G)$ 
of isomorphism classes of rational tropical curves
of toric degree $\Delta$ and combinatorial type $G$ (allowing edge lengths $0$ for convenience)
is parametrised by
\[
	\MM(\Delta, G) \cong \R^n \times (\R_\geq)^{k-3}. 
\]
Here, the factor $\R^n$ parametrizes the position of the marked vertex $h(v_0)$,
and the second factor encodes the lengths of the non-leaf edges of $G$. Again,
there is one exception, namely
\[
  \MM((\delta,-\delta), - \bullet -) \cong \R^n/\R\delta.
\]

\begin{lemma} \label{lem:PointConvergenceHausdorffConvergence}
  Let $h_m : \Gamma_m \to \R^n$ be a sequence of rational tropical curves
	of toric degree $\Delta$ and combinatorial type $G$ converging in 
	$\MM(\Delta, G)$ to a tropical curve $h : \Gamma \to \R^n$.
	Then the sets $h_m(\Gamma_m)$ converge to the Hausdorff limit $h(\Gamma)$.  
\end{lemma}

\begin{proof}
  This follows immediately from the fact that the positions of all vertices 
	and edges of $h_m(\Gamma_m)$ depend linearly (hence continuously) on
	the parameters in $\R^n \times (\R_\geq)^{k-3}$.
\end{proof}

\noindent
Consider the following construction. 
\begin{enumerate}
	\item Mark some of the edges of $G$, including all leaves, by the symbol $\infty$. 
	\item Insert a two-valent vertex in some of the $\infty$-marked edges. 
	      If so, mark both new edges by $\infty$ again.
	\item Decompose $G$ into pieces $G_1, \dots, G_s$ by cutting each interior $\infty$-marked edge
	      into two halves. The ends of the pieces $G_i$ can by canonically labelled by toric 
				degrees $\Delta_i$. 
	\item Mark a vertex $v_i \in G_i$ for $i = 1, \dots, s$. 
	\item Pick a set $S \subset \{1, \dots, s\}$ such that 
	      there exists an assignment of non-negative non-all-zero numbers 
				$(a(e) : e \text{ $\infty$-marked non-leaf})$
				such that for $i,j \in S$ we have 
				\[
				  \sum_{\substack{e \subset [v_i, v_j] \\ \text{$\infty$-marked}}} a(e) \delta(e) = 0.
				\]
				Here, $[v_i,v_j]$ denotes the oriented simple path from $v_i$ to $v_j$.
\end{enumerate}
We call such a construction (and the result $((\Delta_i,G_i)_{i\in S})$) a \emph{degeneration} of $(\Delta, G)$. 
Clearly, $(\Delta_i)_{i\in S}$ is a degeneration of $\Delta$ in the sense of the definition given
before \autoref{thm:TropicalLimits} (set $D$ to be the contraction of $G$ along all non-$\infty$-edges). 
There is an associated linear map of parameter spaces (defined over $\Z$)
\begin{align} %\label{eq:} %\nonumber
  L : \MM(\Delta, G) &\to \prod_{i \in S} \MM(\Delta_i, G_i), \\
	    (\Gamma,h)     &\mapsto (\Gamma_i,h_i)_{i\in S},
\end{align}
given by $h_i(v_i) = h(v_i)$ (identifying $v_i \in G_i$ with $v_i \in G$) 
and keeping the edge lengths
for all edges which are still present. Clearly, the definition
extends in the obvious way to the case when some $G_i$ are $- \bullet - $.
The image $L(\MM(\Delta, G))$ is a rational subcone of $\R^N \times (\R_\geq)^M$ (for suitable $N,M$)
of dimension less than or equal to $\dim(\MM(\Delta, G)) \leq n + k - 3$. 
Moreover, assuming there exist $\infty$-marked non-leaves, the vector $(a(e) : e \text{ $\infty$-marked non-leaf})$
gives rise to a non-trivial kernel element for $L$, hence $\dim L(\MM(\Delta, G)) < \MM(\Delta, G)$.
We can summarize the discussion so far by concluding that in order to prove \autoref{thm:TropicalLimits},
using \autoref{cor:TropicalComplexLimits} it suffices to show the following statement. 

\begin{theorem} %\label{thm }
  Any sequence of tropical rational curves $h_m : \Gamma_m \to \R^n$ of toric degree $\Delta$
	contains a subsequence converging to a Hausdorff limit $A$ (including $A = \emptyset$).  
	The limit $A$ is of the form
	\[
	  A = h_1(\Gamma_1) \cup \dots \cup h_s(\Gamma_s)
	\]
	for a combinatorial type $G$, a degeneration	$((\Delta_1,G_1), \dots, (\Delta_s,G_s))$ of $(\Delta,G)$
	and a tuple of tropical rational curves $(h_1, \dots, h_s) \in L(\MM(\Delta, G))$. 
\end{theorem}

\begin{proof}
  Since the number of trees with $k$ labelled leaves is finite, 
	we can assume that $(h_m)$ has constant combinatorial	type $G$. 
	Throughout the following, we will identify the vertices and edges
	of $G$ with the corresponding vertices and edges of $\Gamma_m$. 
	In particular, given an non-leaf edge $e$ or vertex $v$ of $G$, we write 
	$l_m(e)$ and $h_m(v)$ for the length and position of the corresponding 
	edge and vertex in $\Gamma_m$, respectively. 
	We denote by $\R^n \cup \{\infty\}$ the one-point compactification of $\R^n$. 
	By compactness, we may assume that 
	\begin{itemize}
		\item for any vertex $v \in G$, $h_m(v)$ converges in $\R^n \cup \{\infty\}$,
		\item for any non-leaf edge $e \subset G$, $l_m(e)$ converges in $[0, +\infty]$.
	\end{itemize}
	
	% TpX Bild einfügen
	\begin{figure}[]
		\centering
		\input{pic/DeformationCurve.TpX}
		\caption{The Hausdorff limit of a sequence of tropical rational curves in $\R^2$. On the right hand side,
	the combinatorial type and its degeneration are depicted. The gray parts are the ones we forget in step (e).}
		\label{DeformationCurve}
	\end{figure}%

	We now describe an explicit degeneration of $(\Delta, G)$ (see \autoref{DeformationCurve}). 
	We mark all leaves and all edges with $\lim l_m(e) = +\infty$ by $\infty$ (step (a)). 
	For any such edge $e$, we insert a two-valent vertex if and only if all adjacent vertices diverge and
	there	exists a sequence $x_m \in e^\circ \subset \in \Gamma_m$ such that $h_m(x_m)$ is bounded (step (b)). 
	Passing to a subsequence, we may assume that $h_m(x_m)$ converges in $\R^n$. 
	For each two-valent vertex $v$ we fix such a sequence and set	$h_m(v) = h_m(x_m)$. 
	Let $G_1, \dots, G_s$ denote the pieces after cutting all interior 
	$\infty$-edges into halves (step(c)). 
	We mark a vertex $v_i \in G_i$ for all $i = 1, \dots, s$ (step(d)). 
	Finally, we set	$S = \{i : \lim h_m(v_i) \in \R^n\}$ (step(e)). 
	In other words, we forget all the pieces $G_i$ for which $\lim h_m(v_i) = \infty$. 
	
	Note that since vertices in the same piece $G_i$ are connected via edges with finite 
	limit length, $S$ does not depend on the choice of marked vertices $v_i$. 
	Setting $((h_{i,m})_{i \in S}) = L(h_m)$, we obtain a sequence of tuples
	of rational tropical curves of toric degrees $((\Delta_i)_{i \in S})$ contained
	in $L(\MM(\Delta, G))$. By construction, the limit $\lim L(h_m)$ in 
	$\prod_{i \in S} \MM(\Delta_i, G_i)$ exists. We denote it by $((h_i)_{i \in S})$.
	Since $L(\MM(\Delta, G))$ is closed, it also lies in $L(\MM(\Delta, G))$.
	
	Let us prove that $S$ is allowable.
	Note that the constructed degeneration is non-trivial if and only if 
	it produces at least one interior $\infty$-edge.
	This, in turn, holds true 
	if and only if the sequence $r_m = \max\{l_m(e) : e \text{ non-leaf}\}$ diverges. 
	Then the sequence $(l_m(e) / r_m : e \text{ non-leaf})$ is bounded.
  Let $(a(e) : e \text{ non-leaf})$ denote an accumulation point. 
	Note that $a(e) = 1 \neq 0$	for at least one $\infty$-edge $e$, 
	and that $a(e) = 0$ for any edge not marked by $\infty$. 
	For any pair $i\neq j \in S$, we have
	\[
	  h_m(v_j) - h_m(v_i) = \sum_{e \subset [v_i, v_j]} l_m(e) \delta(e).
	\]
	Dividing by $r_m$ and taking limits, we obtain
	$\sum_{e \subset [v_i, v_j]} a(e) \delta(e) = 0$, as required.
	
	To finish the proof, it remains to show that the sets $h_m(\Gamma_m)$ converge 
	in the Hausdorff sense to 
	\[
	  A = \bigcup_{i \in S} h_i(\Gamma_i).
	\]
	By \autoref{lem:PointConvergenceHausdorffConvergence}, we have
	$\lim_{m \to \infty} h_{i,m} (\Gamma_{i,m}) = h_i(\Gamma_i)$ in the Hausdorff sense. 
	It follows that $\lim A_m = A$ with $A_m = \bigcup_{i \in S} h_i(\Gamma_i)$.
	Let $K \subset \R^n$ be a compact set. For any vertex $v \in G$ which is forgotten 
	during the degeneration construction, we have $h_m(v) = \infty$ and hence
	$h_m(v) \notin K$ for sufficiently large $m$. Let $e \subset G$ be an edge
	which is forgotten during the degeneration. Then $e$ is not subdivided in step (b)
	and both vertices of $e$ converge to $\infty$. If $\lim l_m(e) \neq +\infty$, this 
	implies $h_m(e) \cap K = \emptyset$ for large $m$ by the vertex argument. 
	If $\lim l_m(e) = +\infty$, 
	the same is true since by assumption that there does not exist a sequence of points 
	$x_m$ on $e$ with bounded $h_m(x_m)$. 
	For any other edge $e$, at least one of the adjacent vertices $v$ (possibly after subdividing $e$ 
	into two edges in step (b))  satisfies $\lim h_m(v) \neq \infty$. Then $v \in G_i$ for some 
	$i \in S$ and $h_m(e) \cap K = h_{i,m}(e) \cap K$ for large $m$. 
	It follows that	$h_m(\Gamma) \cap K = A_m \cap K$ for sufficiently large $m$, and the claim follows.
\end{proof}

\printbibliography

\vfill

\section*{Contact}

\begin{itemize}
  \item
    Grigory Mikhalkin, Section de Mathématiques, Université de Genève, Battelle Villa,
	  1227 Carouge, Suisse; \href{mailto:grigory.mikhalkin@unige.ch}{grigory.mikhalkin AT unige.ch}.
  \item
    Johannes Rau, Departamento de Matemáticas, Universidad de los Andes, 
		KR 1 No 18 A-10, BL H, Bogotá, Colombia; \href{mailto:j.rau@uniandes.edu.co}{j.rau AT uniandes.edu.co}.
\end{itemize}

\end {document}